\documentclass[12pt]{article}
\usepackage[english]{babel}
\usepackage{mathptmx}
\usepackage[utf8]{inputenc}
\usepackage{amsmath}
\usepackage{amsfonts}
\usepackage{amssymb}
\usepackage{amsthm}
\newtheorem{thm}{Theorem}[section]

\newtheorem{defn}{Definition}[section]
\newtheorem{lemma}{Lemma}[section]

\usepackage[colorlinks=true,allcolors=blue]{hyperref}
\usepackage{enumerate}
\def\bq{\mathbf{Q}}
\def\bz{\mathbf{Z}}
\title{Exceptional Quartics are Ubiquitous}
\author{Aruna C and P Vanchinathan*}
\begin{document}
\maketitle
\begin{center}
Division of Mathematics\\
VIT University\\
Vandalur–Kelambakkam Road\\
	Chennai, 600 127 INDIA\\[3pt]
	\texttt{
	aruna.2020@vitstudent.ac.in,\quad 
		vanchinathan.p@vit.ac.in}

\end{center}
\begin{abstract}
For each real quadratic field we constructively show the existence of 
infinitely many exceptional quartic number fields containing that
quadratic field.
On the other hand, another infinite collection  of quartic 
exceptional fields without any quadratic subfields is also provided.  
Both these families are non-Galois extensions of $\bq$, and  their  normal closures have Galois groups  $D_4$ and  $S_4$ respectively.
We also show that an infinite number of these exceptional quartic fields
have power integral basis, i.e., monogenic.
We also construct large collections  of exceptional number fields
	in all degrees${}>4$.
\end{abstract}


\section{Introduction}
 Units in the ring of algebraic integers of a number field 
 have been an object of study since Dirichlet.
 A unit $\lambda$ 
 of a number field is said to be exceptional if $1-\lambda$ is also a unit.
  These \textit{exceptional units} were studied extensively by Nagell in a
  sequence of papers \cite{nagell28}, \cite{nagell64}, \cite{nagell68}, \cite{nagell69A}, \cite{nagell69B}.

  A number field that possesses  an exceptional unit is called an
  exeptional number field. Despites this choice of name we want to show
  that exceptional number fields are not really rare to find.
  In a major advancement in 1977, H W Lenstra \cite{lens}
 showed that  a field with sufficient number 
 of exceptional units   such that their differences are
  also units,   will be a norm-Euclidean
  number field.  (See also Houriet \cite{hour}).
  That sufficient number is called the Lenstra constant by later authors
  and following Lenstra they   have found many new number fields that 
  are norm-Euclidean.
  So finding new ways  of constructing   exceptional number fields and 
   and finding  if a large number of them exist are goals worthy of pursuit.
   Both these questions are addressed here with what we believe
   are satisfactory answers.
 
 Before delving into the details of our own work, we list 
 some of the existing major
 results about exceptional number fields.
  \begin{itemize}
	  \item
		  $\bq$($\sqrt{5})$ and $\bq$($\sqrt{-3})$ 
		  are the only exceptional ones among  quadratic fields.
		  (Nagell, 1969 \cite{nagell69A}).
	  \item In a number field there can be at  most a finite number of exceptional units. (Siegel, 1921).
	  \item A number field of degree $n$ has at most $3. 7^{n+2r+2}$ exceptional units, $r$ being the rank of the unit group. (Evertse, 1983.  \cite{ever}).
	  \item General unit equations and 
		  generalization to exceptional 
		  $S$-units and bounds on
		  their heights have been obtained
		  by Gy\H ory \cite{gyory}.
		  (Two whole  monographs 
by Evertse and Gy\H ory and 
		  \cite{gyorybook1}, \cite{gyorybook2}  are devoted to this).
	   \item The only complex cubic fields that are exceptional are
		   the ones with discriminant $-23$ and $-31$.  (\textit{Nagell}).
	  \item The only  exceptional real cubic fields are the following two
		  infinite families, generated by the exceptional 
		  units with  minimal polynomials:  (\textit{Nagell}).
		  \begin{description}
			  \item[Non-Galois:]
    $f_k(x) =x^3+(k-1)x^2-kx-1,\quad    k\in\mathbf{Z},\ k\geq 3$.
   \item[Galois:] $ g_k(x) =x^3+kx^2-(k+3)x+1,\quad    k\in \mathbf{Z},\ k\geq -1$.
\end{description}
	  \item An infinite  family of exceptional quartic number fields, all of unit rank 2,
		  arise from the following irreducible polynomials (Niklasch and Smart
		  \cite{ns}):  
		  \[ x^4+ax^3+x^2+ax-1,\quad a\in\mathbf{Z},\  a>0\]
	  \item Number fields  whose degree is \textit{not} a multiple of 3, where 
		  3 splits  completely
		   do not have exceptional units. (Triantafillou (2021) \cite{tri} ).
	   \item  Among Galois extensions of  $\bq$ 
		   of a given prime degree $p\geq5$,
		   at most finite number of them can be exceptional.
		   In fact, the maximal real subfield
		   of the 11th cyclotomic field is the only exceptional 
		   number field that is cyclic of degree 5.
		   (\cite{frei}

	   \item Leutbecher and Martinet \cite{mart} have constructed
		   many infinite collections of exceptional number fields
		   in degrees 4,5 and 6 and focussed on determining
		   their Lenstra constant so as to find new Euclidean
		   number fields.
  \end{itemize}
  The last mentioned work of 1982 provides the seed for our paper.
We have made a more detailed study of the fields generated from
the two of the 
infinite families of quartic polynomials discovered in \cite{mart}, and also
have constructed many new families of exceptional number fields in degree 4 and above.
(\textit{In fact, unaware of their work we had submitted an earlier version
of our work where
the anonymous referee pointed out this. We thank that referee for providing the reference and also for supplying one argument used in Step (2) of the proof
of part (iv) of Theorem~\ref{main1}}).

Now we summarise our findings below:
\begin{itemize}
	\item [(i)]
  For every \textit{real} quadratic field, there exist
	infinitely many  exceptional quartic number fields containing it.
	\item[(ii)] 
 In these  quartic exceptional 
	fields $\bq(\alpha)$, both  $\alpha$ and $\alpha^2$ are 	exceptional units.
\item[(iii)] Among these exceptional  quartic fields we show that
	infinite subcollection are monogenic, that is, their
ring of integers admit power integral basis.
\item[(iv)] The Galois closures of these quartics  have $D_4$ as Galois groups.
\item[(v)] There is an infinite family of exceptional quartic fields without
	any quadratic subfields.
\item[(vi)] Large families of irreducible polynomials in all degrees$>4$ 
	whose roots are exceptional units are constructed.
\end{itemize}

\subsubsection*{Organization of this paper}
In Section 2 we have put together, with references, all
the preliminary results  which are applied in our paper. 
In Section 3, the quartics  of Leutbecher and Martinet 
are studied in detail proving our main theorem.
 In Section 4, a new  quartic family without quadratic subfields,
 along with higher degree exceptional fields is  constructed.
In the final Section 5, we provide  
exceptional number fields of degree above 4 and mention the avenues
it opens for further research.

 A minor result used in proving our theorem is relegated 
 to an appendix so as not to
hinder the flow.

\section{Preliminaries}
In this section we assemble  various definitions, and results of earlier 
authors -- some of them more than a century old -- that are needed to prove our main theorem.
\begin{defn}
Let $F$  be a number field. A unit $\lambda$  of the ring of integers of 
$F$ is said to be an exceptional unit when $1- \lambda$ is also a unit. 
A number field is said to be exceptional if its ring of integers 
contains an exceptional unit.
\end{defn}

 The following easy-to-verify result provides us with a simple  condition for 
 a polynomial to have exceptional units as its roots.  
\begin{lemma}\label{excep}(Nagell, 1969)
Let $f(x)\in \mathbf{Z}[x]$ be a monic irreducible polynomial. If $|f(0)|=1$ and $|f(1)|=1$, then any root $\alpha$ of $f(x)$ is an exceptional unit in $\bq(\alpha)$.
\end{lemma}
Besides the well-known arithmetic condition of Eisenstein for the irreducibility of polynomials there are some other non-arithmetic criteria. The following sufficient condition for irreducibility by Perron will be used 
here in many of  our  constructions.
\begin{thm}[Perron]\label{perron}
Let
\begin{equation*}
    f(x)=x^n+a_1x^{n-1}+\cdots+a_n
\end{equation*}
 be a polynomial with integer coefficients such that $a_n \neq0$. Further
 \begin{enumerate}[(i)]
     \item if $|a_1| > 1+|a_2|+\cdots+|a_n|$,  then f is irreducible, or
     \item if $|a_1| \geq 1+|a_2|+\cdots+|a_n|$ and $f(\pm1)\neq0$,
	     then f is irreducible.
 \end{enumerate}
\end{thm}
Proof of this result makes  use of complex analysis, and  can be found in the monograph of
	 V.V. Prasolov\cite{pra}.

The following classical result\cite{rees} will be used later for proving our theorem.
\begin{lemma}\label{realroot}
Let $f(x)=x^4+ax^3+bx^2+cx+1$ be a monic polynomial with real coefficients. The discriminant and two other useful invariants for this polynomial are below:
\begin{align*}
    \Delta &=256-192ac-128b^2+144bc^2-27c^4+144a^2b-6a^2c^2-80ab^2c\\
	&{} \hspace{2.83em}{}   +18abc^3+16b^4-4b^3c^2-27a^4+18a^3bc-4a^3c^3-4a^2b^3+a^2b^2c^2\\
    D & =64-16b^2+16a^2b-16ac-3a^4\\
    P & =8b-3a^2 
\end{align*}
	Then a sufficient condition on $f(x)$ to have all its roots real and distinct is that
 $\Delta>0$, $P<0$, and $D<0$.

\end{lemma}

\begin{defn}
    Let $f(x;t)$  be a monic irreducible polynomial in $\mathbf{Z}[x]$. Let $\Delta(t)\in \mathbf{Z}[t]$ be the discriminant of f with respect to x. We have $\Delta(t)=\prod_{i=1}^{r}p_i(t)^{e_i}$ with pairwise distinct irreducible elements $p_i(t)$ in $\mathbf{Z}[t]$. Then $\Delta^{\rm red}(t):=\prod_{i=1}^{r}p_i(t)$ is called the reduced discriminant of $f$.
\end{defn}

The following result of J.K\"onig with a flavour of 
Hilbert's Irreducibility Theorem will be used:
\begin{lemma}\label{monogenic}(See \cite{joachim})
    Let $h(t,x):=h_1(x)-th_2(x)\in \mathbf{Z}[t,x]$ be a polynomial with Galois group $G$ over $\bq(t)$, and $\Delta^{\rm red}(t)\in \mathbf{Z}[t]$ be its reduced discriminant. Assume that the following hold:
    \begin{enumerate}[(i)]
	    \item $\Delta^{\rm red}(t)$ has no irreducible factor of degree${}\geq4$.
        \item There is no prime $p$ that divides all the values in $\Delta^{\rm red}(\mathbf{Z})$.
    \end{enumerate}
    Then, there are infinitely many choices for  $t_0 \in \mathbf{Z}$ such that
	choosing  $\alpha$ as a root of $h(t_0,x)=0$, gives 
	a monogenic number field $\bq(\alpha)$  with Galois group
	isomorphic to $G$.
More precisely, $\{\alpha^i\mid i= 0,1,\ldots,\deg(h)-1\}$ is a  power basis of the ring of integers of 
	$\bq(\alpha)$.
\end{lemma}

\section{Two Quartic Families with Quadratic Subfields}
With all the preliminaries in place we are now ready to state and prove the
first one of our main results.

First we look at the family  polynomials denoted by
$P_a(x)$ by Leutbecher and Martinet \cite{mart} and 
studied to determine their Lenstra constant. Our theorem below 
shows they have many remarkable properties.
\begin{thm}\label{main1}
For an integer $t\geq4$ , consider the following polynomial
\begin{equation}
  f(x;t)=x^4-tx^3-x^2+tx+1
\end{equation}
Then we have
\begin{enumerate}[(i)]
    \item $f(x;t)$ is irreducible 
    with all its roots  real, thus defining a unit rank 3 
number field $K_t=\bq(\alpha)$ where $\alpha$ is a root of $f(x;t)$.
    \item $K_t$ is an exceptional number field where both 
	    $\alpha$ and $\alpha^2$ are exceptional units.
    \item   $K_t$ contains  the real quadratic field $\bq(\sqrt{t^2-4})$.
    \item Every real quadratic field $\bq(\sqrt d)$ is contained as a subfield
	    in  infinitely many  quartic fields $K_t$.
    \item The  Galois group of $f(x;t)$ is 
	    the dihedral group $D_4$.
\end{enumerate} 
\end{thm}
Proof of (i): 
Obviously  $f(x;t)$ has no rational root. So if it were reducible then 
it ought to be a product of two degree 2 polynomials
which can be ruled out by routine calculations.
About its roots, we will appeal to Lemma~\ref{realroot} and use the 
notations in that lemma.
\begin{align*}
    \Delta & =4t^6-23t^4-8t^2+144 =(t^2-4)^2(4t^2+9)>0\\
    P & =-8-3t^2 <0\\
	D & =-3t^4+48<0\quad \mbox{ for }  t\geq 3
\end{align*}

Clearly $\Delta>0,\ P<0$ and $D<0$ for all $t\in\mathbf{Z}\geq 3$. So by Lemma \ref{realroot}, all the roots of  the polynomial 
$f(x;t)$ are real. Hence the unit rank of the number field generated by 
any root of the polynomial is 3.


\vspace{6pt}
Proof of (ii):
    We have $f(0)=1$ and $f(1)=-1$, so by Lemma \ref{excep}, the root of this polynomial is an exceptional unit. Hence the number field generated by any
    root $\alpha$  of the polynomial is an exceptional number field.
    We have now $1-\alpha$ is a unit.

    One can easily verify $1+\alpha$ is also a  unit by directly 
    calculating its minimal polynomial.
    Now the product of the two units
    $1-\alpha$ and $1+\alpha$, that is,  $1-\alpha^2$,  has to be a unit,
    thereby proving  $\alpha^2$ is an exceptional unit.
    
\vspace{6pt}
   Proof of (iii):   From the minimal polynomial satisfied by $\alpha$,
   by simple rearrangement we get
 $$\alpha^4-\alpha^2 +1 = t (\alpha^3-\alpha),$$
which in turn gives us
   $$t=\frac{\alpha^2-1}{\alpha}+\frac{\alpha}{\alpha^2-1}.$$
 Now denoting the quantity $\frac{\alpha^2-1}{\alpha}$ by 
 $\beta $ (which is in $ \bq(\alpha))$, 
   we get $$t=\beta+\frac{1}{\beta}.$$
   From this we  $t^2-4=(\beta-\frac{1}{\beta})^2$  showing 
   the integer $t^2-4$ to be   a square in the
   quartic field $K_t=\bq(\alpha)$.
   As $t^2-4$ is never a perfect square in $\bz$, we see that 
   the field $\bq(\sqrt{t^2-4}) $ indeed is a quadratic subfield of $K_t$.

\vspace{6pt}
   Proof of (iv): Step 1. Initially we just show for each real quadratic 
there is one $t$ such that $K_t$ has that field as a subfield.
   
   Given a positive square-free integer $d$,
   and a unit $u$ in  the  real quadratic field 
   $\bq(\sqrt d )$ writing $2u= t+s\sqrt d$, for some $t,s\in\bz$
   and computing its norm we get  $t^2-ds^2 =4$, that is,    $t^2-4 = ds^2$.
   The last equality shows 
   that the two numbers,  $\sqrt d$ and $\sqrt{t^2-4}$,
   are merely two different primitive elements
   for the same  real quadratic  extension of $\bq$.
   This, combined with part (iii) just proved above, shows that 
   $\bq(\sqrt d)$ is a subfield of  our  exceptional quartic
   $K_t$. Let us denote this $t$ as $t_1$.

   Step 2: Now we define an infinite sequence of positive 
   integers $t_i$, $i=2,3,\ldots$
   such that all  these $K_{t_i}$ contain the same quadratic field
   as $K_{t_1}$ does.
   For this purpose we  define inductively,
   $t_{i+1} = t_i^2-2$.
   The quadratic subfield  of $K_{t_{i+1}}$ by (iii), is
   the one generated by $\sqrt{t_{i+1}^2-4}$ which by
   definition of $t_{i+1}$ is $\sqrt{ (t_i^2-2)^2-4} =
   \sqrt{ t_i^2(t_i^2-4) }$ which is just another primitive element 
   for the quadratic subfield of $K_{t_i}$.

   Though we have shown infinitely $t_i\in\bz$ exists such that
   $\bq(\sqrt d)\subset K_{t_i}$ we still need to show that
   infinitely many of these  $K_{t_i}$ are distinct fields.

   We already know that the roots of the irreducible 
   polynomials $f(x,t_i)$ are exceptional units in $K_{t_i}$, and they
   are distinct for  two different values of $t_i$, being roots
   of two different irreducuble polynomials. Now combining
   with the fact that there are only finitely many exceptional
   units in a single number field, we can then conclude that
   infinitely many 
   among the  quartic  fields $K_{t_i}$ are distinct, completing the
   proof of part (iv).

Proof of (v): For this we will make use of the well-known
results on determination of Galois groups of quartics
 found, for example, in \cite{kappe}. 
Accordingly we will calculate the discriminant and the resolvent for $f(x,t)$.
    The discriminant 
    is $\Delta=(4t^2+9)(t^2-4)^2$. Clearly, $\Delta$ can be a square only when $4t^2+9$ is a square. It is easy to see that $4t^2+9$ is not a square for $t\geq3$.
    \par The resolvent polynomial for our $f(x;t)$, denoted $R_3(x)$ in
    \cite{kappe}
  is easily calculated to  be  $R_3(x)=x^3+x^2-(t^2+4)x-(2t^2+4)$,
 and   factorizes as $(x+2)(x^2-x-t^2-2)$.
 So the Galois group is $D_4$ or $\mathbf{Z}/4\mathbf{Z}$.
    \par Note that for $R_3(x)$,  $-2$ is  a rational root.
    So continuing to follow \cite{kappe}, 
 we will be able to conclude  that the Galois group is
$D_4$ once we  show that 
   $(t^2-4)\Delta$  is not a square.
This is equivalent to showing $(t^2-4)(4t^2+9)$ is not a square, for
$t\geq3$.
In the appendix we have included a proof of this.
  Hence the Galois group of the splitting field of $f(x;t)$ is $D_4$.

\vspace{6pt}
Next we move onto  discuss second  family of irreducible
quartic polynomials denoted by $R_a(x)$ in \cite{mart}.
Their  roots are exceptional units.
Compared with the first  family the difference is only in
the coefficient of $x^2$. Here are the results of our detailed study
of the number fields corresponding to them:
\begin{thm}
For an integer $t\geq$ 7, the following polynomial
\begin{equation}
  h(x;t)=x^4-tx^3-3x^2+tx+1
\end{equation}
is irreducible.
	Any root $\alpha$ of $h(x,t)$ is exceptional,
	and $1+\alpha$ is also a unit, making $\alpha^2$ exceptional.
	The exceptional quartic number fields generated
	by $\alpha$ contain the quadratic subfield
	$\bq(\sqrt {t^2+4})$, and their normal closure
	have $D_4$ as Galois groups.
The number field $Q(\alpha)$ is monogenic for infinitely many integers t.
\end{thm}
We will prove only the monogenity part as other
proofs run along the same lines as in the first family. For this
we appeal to the result of K\"onig, quoted as Lemma~\ref{monogenic} 
in Section 2.

The discriminant of the polynomial is $\Delta=(4t^2+25)(t^2+4)^2$,
so (in the notation of  Lemma~\ref{monogenic}),
    $\Delta^{\rm red}(t)=(4t^2+25)(t^2+4)$ 
  which has no irreducible factor of degree $\geq4$.

  Further the value for the reduced discriminant at 
  $t=0$   and $t=3$ are 100 and 793 respectively, 
  two numbers with  no common factors. 
  So there is no prime
    $p$ dividing $\Delta^{\rm red}(\bz)$. Hence there are infinitely many integers $t$ such that $\bq(\alpha)$ is a monogenic number field by Lemma \ref{monogenic}. 

 \noindent\textbf{Remark:}\quad  Given an $\alpha$, an exceptional unit, we  have some more naturally associated exceptional units, namely $\frac{1}{\alpha}, 1-\alpha, \frac{1}{1-\alpha}, 1-\frac{1}{\alpha}, 1-\frac{1}{1-\alpha}$.
 Also in both the polynomials $f(x,t), h(x,t)$ whenever $\alpha$ is a root
 $-1/\alpha$ is also a  root, which is again an exceptional unit. 
 Along with $\alpha^2$, this gives rise
 to 18 exceptional units  in all the fields from either  families.
 Also, notice that both $\alpha,\alpha^2$ being  exceptional units, allows
 us to conclude that $1,\alpha,\alpha^2$ can be part of an exceptional
 sequence in the terminlogy of Lenstra.


 \section{Higher Degree  Family where the  Quartics 
 have no Subfields}

In the previous two families the coefficients of $x^2$ was independent
of the parameter $t$ and that helped us prove their
their irreducibility by 
direct ad-hoc calculations. This time no coefficients, with obvious
exceptions, will be fixed  and we will
appeal to Perron's irreducibility theorem.
\begin{thm}\label{main2}
For $n,t\geq4$ both  integers
	consider
the collection of polynomials:
\begin{equation}
    g_n(x;t)=x^n-(t+3)x^{n-1}+tx+1
\end{equation} 
Then we have 
 \begin{enumerate}[(i)]
     \item each $g_n(x;t)$ is irreducible;
     \item any root of $g_n(x;t)$ defines an exceptional number field;
     \item for $n=4$, the polynomial $g_4(x;t)$ defines a quartic number field $L$ of unit rank 3 and the polynomial has  Galois group  $S_4$;
     \item the above quartic number field $L$
		 has no quadratic subfield.
 \end{enumerate}
\end{thm}
\vskip2pt
\begin{proof}[Proof of (i)]
To apply Perron's theorem \ref{perron} we need to check an inequality
involving the coefficients. Here we need to check 
	 $|t+3|>1+t+1$ which is obviously true showing
	$g_n(x)$ is irreducible.
\end{proof}
\begin{proof}[Proof of (ii)]
	Clearly $g_n(0)=1$ and $g_n(1)=-1$, 
and so by Lemma \ref{excep} 
	 any root of $g_n(x)$ generates an exceptional number field.
\end{proof}
\begin{proof}[Proof of (iii)]
To prove that all roots are real we again appeal to Lemma \ref{realroot} and use notations in that lemma.
\begin{align*}
    \Delta & =4t^6+36t^5+48t^4-252t^3-1320t^2-2340t-1931 > 0\\
    P & =-3(t+3)^2 <0\\
    D & =-3(t+3)^4+16t+112<0
\end{align*}
Clearly $P<0$ and $D<0$ for any $t$ but $\Delta<0$ only for $t\geq4$. This implies $g_4(x;t)$ has all its roots real for $t\geq4$. So the unit rank of the number field generated by the root of the polynomial is 3.

\par Now we move on to prove the Galois group of $g_4(x;t)$ is $S_4$.
	Again we appeal to the criterion given in \cite{kappe}.
	In this case
	  the resolvent polynomial is easily seen to be 
	$R_3(x)=x^3-(t^2+3t+4)x-(2t^2+6t+9)$. 
	Routine calculations show that this is irreducible.
	So the Galois group has to be either $S_4$ or $A_4$.\\

	The polynomial $g_4(x)$ is easily checked to be irreducible mod 2, 
	showing the existence of a 4-cycle in the Galois group.
	So the Galois group is $S_4$. 
\end{proof}

Proof of (iv): 
    Note that  $S_3$ is a maximal subgroup in $S_4$. Hence, by the fundamental Galois correspondence, the quartic number field being the fixed 
    field of $S_3$ in an $S_4$-Galois extension,  cannot have a 
    quadratic subfield.
\section{More Higher Degree Exceptional fields}
 
 Perron's criterion for irreducibility makes it easy
to produce irreducible polynomials whose roots will be
 exceptional units of any desired degree. 

The problem is how to infer anything about the structure of these
fields just from this polynomial.
 For example the roots of the
following polynomials are all exceptional units:
\begin{itemize}
	\item $x^n-(t+3)x^{n-1}+tx^k+1\qquad  1\le k\le n-2$
	\item $x^n-(t+n-2)x^{n-1}+tx^{n-2}+x^{n-3}+\cdots+x+1$
\end{itemize}
Another  family of polynomials 
in multiple parameters 
$F(x;t_1,t_2,\ldots, t_{n-2})$
of degree $n$ for a given choice of $n-2$
ordered  integer coefficients,  $t_1,  t_2,\ldots, t_{n-2}$ 
is below:
\begin{equation*}
F(x;t_1,t_2,\ldots, t_{n-2})
	=x^n-(t_1+t_2+\cdots+t_{n-2}+3)x^{n-1}+t_1x^{n-2}+t_2x^{n-3}+\cdots+t_{n-2}x+1
\end{equation*}
It is intriguing to note  that any permutation of these $n-2$ coefficients 
produces irreducible  polynomials with all their roots  exceptional units of degree $n$.
It remains to be investigated  what nice properties these exceptional number 
fields will have when the coefficients $t_i$'s, and their ordering,
are chosen carefully.

\section*{Appendix}
Now we provide a proof of the statement that was used to prove
certain polynomial has $D_4$ is its Galois group.

\noindent\textit{For $t\in\bz$, if $(t^2-4)(4t^2+9)$ is a perfect square then
$t=3$.}

Proof:\quad First note that $t^2-4$ itself is never a square, so one can
write it uniquely as $t^2-4 = dn^2$ where $d$ is the square-free part of
it, $d>1$, and $n\in \bz$. 
First we rewrite the number under discussion as below:
 $$(t^2-4)(4t^2+9)= (t^2-4)[ 4(t^2-4) +25\big].$$
Let us assume this number is a square $N^2$.
And substituting for $t^2-4$ gives us
$$dn^2(4dn^2+25)=N^2\ldots\ldots\ldots(1)$$
As $d$ is square-free, and the rhs is a square it follows that
$d$ is a divisor of the quantity inside the  bracket, namely $4dn^2+25$,
and hence $d|25$. But $d$ being  square-free and $>1$, 
 the only conclusion is $d=5$.

So Eq.(1) becomes
$5n^2(20n^2+25)=N^2$, and hence $25n^2(4n^2+5)=N^2$.
As $25n^2= (5n)^2$, a perfect square, 
we can conclude  that $4n^2+5$ is  a square
which is possible only when $n=1$.
That makes $t^2-4=dn^2 = 5.1^2 = 5$, and so $t=3$ is the
only solution. \hfill QED



\frenchspacing
\begin{thebibliography}{99}
	\bibitem{frei}N. Freitas, A.Kraus, S.Siksek,
		On the unit equation of cyclic number fields of prime degree,
		\textit{Algebra and Number Theory}, Vol. 15, (2021),
		  2647--2653.
	\bibitem{gaal}I. Ga\'al, \textit{Diophantine Equations and Power Integral Bases}, Birkhauser, 2019.
\bibitem{gyory} K. Gy\H ory, Bounds for the solutions of $S$-unit equations and
	decomposable form equations II, \textit{Publ. Math. Debrecen}, 2019,   507--526.
\bibitem{gyorybook1}J.-H. Evertse and K.Gy\H ory,\textit{Unit Equations in Diophantine Number Theory}, Cambridge Studies in Advanced Mathematics Series, Vol. 146, Cambridge, 2015. 
\bibitem{gyorybook2}J.-H. Evertse, K. G\H yory,
	\textit{Discriminant Equations in Diophantine Number Theory}. New
Mathematical Monographs, vol. 32 (Cambridge University Press, Cambridge, 2017).
\bibitem{hour} J. Houriet, Exceptional units and Euclidean number fields, \textit{Arch.Math.} 88, 2007,  425--433.
\bibitem{ever}J.-H. Evertse, Upper bounds for the number of solutions of Diophantine equations, \textit{Amterdam Stichting Mathematisch Centrum}, 1983.
\bibitem{kappe}  L.C.Kappe and  B.Warren,  An elementary test for the Galois group of a quartic polynomial, \textit{American Mathematical Monthly}, 96, 1989,   133--137.
\bibitem{joachim} J. K\"onig, A note on families of monogenic number fields.
\textit{Kodai J. Math}, 41, 2018,   456--464.
\bibitem{mart}A. Leutbecher and J. Martinet, Lenstra's Constants and Euclidean
	Number Fields, \textit{Ast\'erisque}, Tome 94, (1982), pp. 87--131.
\bibitem{lens} H.W. Lenstra, Euclidean number fields of large degree, \textit{Invent. Math.} 38,   237--254.
\bibitem{nagell28} T. Nagell, Darstellung ganzer Zahlen durch bin\"are kubische Formen mit negativer Diskriminante, \textit{Mathematische Zeitschrift}, 28, 1928,   10--29.
\bibitem{nagell64} T. Nagell, Sur une propri\'et\'e des unit\'es d’un corps alg\'ebrique, \textit{Arkiv f$\ddot{o}$r Matematik}, Vol. 5, 1964,   343--356.
\bibitem{nagell68} T. Nagell, Sur les unit\' es dans les corps biquadratiques primitifs du premier rang,\textit{Arkiv f$\ddot{o}$r Matematik}, Vol. 7, 1968,   359--394.
\bibitem{nagell69A} T. Nagell, Quelques probl\` emes relatifs aux unit\'es alg\'ebriques, \textit{Arkiv f\"or Matematik} 8 (1969),   115--127.
\bibitem{nagell69B}  T. Nagell, Sur un type particulier d’unit\' es algébriques, \textit{Arkiv f$\ddot{o}$r Matematik} 8, 1969,   163--184.
\bibitem{nick} G. Niklasch,  Family portraits of exceptional units,
\bibitem{ns} G. Niklasch,  and  N.P. Smart, Exceptional Units in a Family of Quartic Number Fields, \textit{Mathematics of Computation}, Vol. 67, No. 222, 1998,   759--772.
\bibitem{pra} V. V. Prasolov, \textit{Polynomials}, 
	Algorithms and computation in mathematics, Vol. 11, Springer
\bibitem{rees} E.L. Rees, Graphical Discussion of the Roots of a Quartic Equation, \textit{The American Mathematical Monthly}, Vol. 29, No. 2, 1922,   51--55.
\bibitem{siegel} C. L. Siegel, \"Uber einige Anwendungen diophantischer Approximationen,
	\textit{Abh. Preuss. Akad.Wiss}, 1929,   1--41.
\bibitem{tri} N.Triantafillou, There are no Exceptional Units in Number Fields of degree prime to 3 where 3 splits completely, \textit{Proceddings of American Mathematical Society}, Series B,  Vol. 8, 2021,   371--376.
\end{thebibliography}
\end{document}